\input amstex
 \input epsf
\magnification=\magstep1 
\baselineskip=13pt
\documentstyle{amsppt}
\vsize=8.7truein \CenteredTagsOnSplits \NoRunningHeads

\topmatter
\title  Approximating real-rooted and stable polynomials, with combinatorial applications \endtitle 
\author Alexander Barvinok  \endauthor
\address Department of Mathematics, University of Michigan, Ann Arbor,
MI 48109-1043, USA \endaddress
\email barvinok$\@$umich.edu  \endemail
\date June  2018 \enddate
\thanks  This research was partially supported by NSF Grant DMS 1361541.
\endthanks 
\abstract Let $p(x)=a_0 + a_1 x + \ldots + a_n x^n$ be a polynomial with all roots real and satisfying $x \leq -\delta$ for some $0<\delta <1$. We show that for any 
$0 < \epsilon <1$, the value of $p(1)$ is determined within relative error $\epsilon$ by the coefficients $a_k$ with $k \leq {c  \over \sqrt{\delta}}  \ln {n \over \epsilon \sqrt{ \delta}}$
for some absolute constant $c > 0$. Consequently, if $m_k(G)$ is the number of matchings with $k$ edges in a graph $G$, then for any $0 < \epsilon < 1$, the total number $M(G)=m_0(G)+m_1(G) + \ldots $ of matchings is determined within relative error $\epsilon$ by 
the numbers $m_k(G)$ with $k \leq c \sqrt{\Delta} \ln (v /\epsilon)$, where $\Delta$ is the largest degree of a vertex, $v$ is the number of vertices of $G$ and $c >0$ is an absolute constant.
We prove a similar result for polynomials with complex roots satisfying $\Re\thinspace z \leq -\delta$ and apply it to estimate the number of unbranched subgraphs of $G$.
\endabstract
\endtopmatter

\document

\head 1. Introduction and main results \endhead

Our main motivation comes from the observation that in some cases, the total number of combinatorial structures of a particular type is determined with high accuracy by the exact number of the structures of the same type but of a small (sometimes, very small) size. For example, if $G$ is a graph and we are interested in the number $M(G)$ of {\it matchings} in $G$, that is, collections of vertex-disjoint edges, then within a given relative error $0 < \epsilon <1$, the number $M(G)$ is determined by the numbers $m_k(G)$ of matchings with exactly $k$ edges for $k \leq c \sqrt{\Delta(G)} \ln {v(G) \over \epsilon}$, where $\Delta(G)$ is the largest degree of a vertex of the graph, $v(G)$ is the number of vertices of $G$ and $c> 0$ is an absolute constant. We deduce this, and some related results, from some general observation about real-rooted polynomials.

Below we talk about approximating some real and complex values up to ``relative error $\epsilon$". Given a complex number $a \ne 0$, we say that a complex number $b \ne 0$ approximates $a$ up to (or within) relative error $\epsilon>0$ if we can write $a=e^z$ and $b=e^w$ for some complex numbers $z$ and $w$ such that $|z-w| \leq \epsilon$.

We prove the following main result.
\proclaim{(1.1) Theorem} Suppose that $p(x)=a_0 +a_1 x + \ldots + a_n x^n$ is a polynomial and all roots $x$ of $p$ are real and satisfy 
$x \leq -\delta$ for some $0 < \delta < 1$. Then, for any $0 < \epsilon < 1$, the value of $p(1)$, up to  relative error $\epsilon$, is determined by the coefficients 
$a_k$ with 
$$k \ \leq \ {c \over \sqrt{\delta}} \ln {n \over \epsilon  \sqrt{\delta}},$$
where $c >0$ is an absolute constant.
\endproclaim 
In fact, we present a polynomial time algorithm, which, given $a_0, \ldots, a_k$ with $k$ as in the theorem, computes $p(1)$ within relative error $\epsilon$. We note that we need to know $a_0, \ldots, a_k$ exactly (or, at the very least, with very high precision).
Theorem 1.1 is a much improved version of the ``personal communication"  of the author that was referred to in \cite{PR17}.

Using the Heilmann-Lieb Theorem \cite{HL72},  see also \cite{GG81}, we immediately deduce the result about matchings in graphs. We consider undirected graphs, without loops or multiple edges.

\proclaim{(1.2) Theorem} For a graph $G$, let $m_k(G)$ be the number of matchings that contain exactly $k$ edges and let $M(G)=m_0(G) + m_1(G)+ \ldots $ be the total number of matchings. Then, for any $0 < \epsilon <1$, up to relative error $\epsilon$, the number $M(G)$ is determined by the numbers 
$m_k(G)$ with 
$$k \ \leq \ c \sqrt{\Delta(G)} \ln {v(G) \over \epsilon},$$
where $\Delta(G)$ is the largest degree of a vertex of the graph, $v(G)$ is the number of vertices of $G$ and $c> 0$ is an absolute constant.
\endproclaim
Again, we have a polynomial time algorithm, which, given $m_k(G)$, with $k$ as in the theorem, produces an estimate of $M(G)$ within relative error $\epsilon$.
We note that Patel and Regts \cite{PR17} constructed a polynomial time algorithm for computing $m_k(G)$ with $k=O(\ln v(G))$ provided the largest degree $\Delta(G)$ is 
fixed in advance (a straightforward enumeration gives only a quasi-polynomial algorithm of $v(G)^{O(\ln v(G))}$ complexity). For general graphs, the complexity of the algorithm roughly matches that of Bayati et al. \cite{B+07}, which estimates $M(G)$ using the correlation decay approach. Although our approach and that of \cite{B+07} look completely different, they are both inspired by the concept of ``phase transition" coming from statistical physics; more precisely, two related, but different concepts: ours has to do with the Lee-Yang approach via complex zeros of the ``partition function" \cite{YL52}, \cite{LY52} while that of \cite{B+07} has to do with correlation decay, cf. \cite{DS87} and \cite{KK98}. Hence the fact that the complexity appears to be roughly the same is not entirely accidental.

Using the Chudnovsky-Seymour extension \cite{CS07} of the Heilmann-Lieb Theorem and the Dobrushin-Shearer  bound on the roots of the independence polynomial, see, for example, \cite{SS05}, we get another combinatorial application.  Recall that a subset of vertices of a graph is called an {\it independent set} if no two vertices of the subset span an edge of the graph. A graph is called {\it claw-free} if it does not contain an induced subgraph consisting of a vertex connected to some other three vertices that are pairwise unconnected. We obtain the following result.

\proclaim{(1.3) Theorem} For a graph $G$, let $i_k(G)$ be the number of independent sets with exactly $k$ vertices and let $I(G)=i_0(G) + i_1(G)+ \ldots $ be the total number of independent sets. Then, for any $0 < \epsilon <1$, up to relative error $\epsilon$, the number $I(G)$ of a claw-free graph is determined by the numbers 
$i_k(G)$ with 
$$k \ \leq \ c \sqrt{\Delta(G)} \ln {v(G) \over \epsilon},$$
where $\Delta(G)$ is the largest degree of a vertex of the graph, $v(G)$ is the number of vertices of $G$ and $c> 0$ is an absolute constant.
\endproclaim
We have a polynomial time algorithm, which, given $i_k(G)$, with $k$ as in the theorem, produces an estimate of $I(G)$ within relative error $\epsilon$. Curiously, while 
the correlation decay approach of \cite{B+07} is essentially harder in the case of independent sets in claw-free graphs than it is in the case of matchings, our approach is the same in both cases (assuming, of course, the hard work done in \cite{HL72} and \cite{CS07}).

Theorems 1.2 and 1.3 vaguely resemble the ``approximate inclusion-exclusion" of \cite{LN90} and \cite{K+96}. The methods, however, look completely different. It would be interesting to find out if there is indeed any connection between our Theorem 1.1 and the results of \cite{LN90} and \cite{K+96}.

Next, we consider polynomials $p(z)=a_0 +a_1 z + \ldots + a_n z^n$ with complex roots satisfying $\Re\thinspace z \leq - \delta$ for some $0 < \delta <1$ (we call such polynomials ``stable"). We allow complex coefficients $a_k$. We obtain the following result.

\proclaim{(1.4) Theorem} Suppose that $p(z)=a_0 +a_1 z + \ldots + a_n z^n$ is a complex polynomial and all roots $z$ of $p$ satisfy 
$\Re\thinspace z \leq -\delta$ for some $0 < \delta < 1$. Then, for any $0 < \epsilon < 1$, the value of $p(1)$, up to  relative error $\epsilon$, is determined by the coefficients 
$a_k$ with 
$$k \ \leq \ {c \over \delta} \ln {n \over \epsilon \delta},$$
where $c >0$ is an absolute constant.
\endproclaim 

We apply Theorem 1.4 to count {\it unbranched subgraphs}, that is, collections of edges of the graph such that every vertex of the graph is incident to at most two edges from the collection. From  Ruelle's Theorem \cite{R99a}, \cite{R99b}, see also \cite{Wa09}, we deduce the following result.

\proclaim{(1.5) Theorem} For a graph $G$, let $u_k(G)$ be the number of unbranched subgraphs with exactly $k$ edges and let $U(G)=u_0(G) + u_1(G)+ \ldots $ be the total number of unbranched subgraphs. Then, for any $0 < \epsilon <1$, up to relative error $\epsilon$, the number $U(G)$ is determined by the numbers 
$u_k(G)$ with 
$$k \ \leq \ c \left(\Delta(G)\right)^3 \ln {v(G) \over \epsilon},$$
where $\Delta(G)$ is the largest degree of a vertex of the graph, $v(G)$ is the number of vertices of $G$ and $c> 0$ is an absolute constant.
\endproclaim
One can easily see that if a non-constant polynomial $p$ satisfies the conditions of Theorems 1.1 or 1.4 then so does its derivative $p'$. Therefore, in Theorems 1.2, 1.3 and 1.5, we can not only estimate the number of structures of a given type (matchings, independent sets or unbranched subgraphs) by counting structures up to some small size, but also estimate the average size of a structure, the second moment, etc.

Finally, we mention the following result implicit in Section 2.2 of \cite{Ba16}. 
\proclaim{(1.6) Theorem} Let $p(z)=a_0+a_1 z + \ldots + a_n z^n$ be a polynomial and suppose that for some $0 < \delta <1$ we have 
$p(z) \ne 0$ for all $z$ in the $\delta$-neighborhood of the interval 
$[0, 1] \subset {\Bbb C}$ (we measure distances by identifying ${\Bbb C}={\Bbb R}^2$). Then, for any $0 < \epsilon <1$, the value of $p(1)$, up to relative error $\epsilon$, is determined by the coefficients $a_k$ with 
$$k \ \leq \   \exp\left\{O\left({1 \over \delta}\right)\right\}  \ln { n \exp\left\{ O\left( {1 \over \delta}\right)\right\} \over \epsilon}.$$
\endproclaim
For applications of Theorem 1.6 to computing partition functions, see \cite{Ba16}. 

We can replace the exponential dependence on $1/\delta$ in Theorem 1.6 by a polynomial dependence if we assume that $p(z) \ne 0$ for $z$ in the $\delta$-neighborhood 
of the sector $|\arg z | < \alpha$ for some fixed $\alpha >0$ and some $\delta >0$. We briefly discuss this in Section 2 and applications to counting subgraphs with prescribed degrees  in Section 3.

We prove Theorems 1.1 and 1.4 in Section 2 and Theorems 1.2, 1.3 and 1.5 in Section 3.

\head 2. Proofs of Theorems 1.1 and 1.4 \endhead

We denote the complex plane by ${\Bbb C}$, the Riemann sphere ${\Bbb C} \cup \{\infty\}$ by $\widehat{\Bbb C}$ and the open unit disc by ${\Bbb D}$, so that 
$${\Bbb D}=\left\{z \in {\Bbb C}:\ |z| < 1 \right\}.$$ By $c$ we denote a positive absolute constant, whose value may change from line to line.

We start with a couple of lemmas.

\proclaim{(2.1) Lemma} Let $h_1, h_2: {\Bbb C} \longrightarrow {\Bbb C}$ be polynomials of degrees $n_1$ and $n_2$ respectively and let 
$$g(z)={h_1(z) \over h_2(z)}, \quad g: \widehat{\Bbb C} \longrightarrow \widehat{\Bbb C},$$
be a rational function. Let $\beta > 1$ be a real number and suppose that 
$$h_1(z) \ne 0 \quad \text{and} \quad  h_2(z) \ne 0 \quad \text{provided} \quad |z| < \beta,$$
so $g$ has neither zeros no poles in the disc $\beta {\Bbb D}=\left\{z \in {\Bbb C}:\ |z| < \beta \right\}$.

Let us choose a branch of
$$f(z)=\ln g(z) \quad \text{where} \quad |z| \ < \ \beta$$
and let 
$$T_m(z)=f(0)+\sum_{k=1}^m {f^{(k)}(0) \over k!} z^k$$
be the Taylor polynomial of degree $m$ of $f(z)$ computed at $z=0$. Then 
$$\left| f(1)-T_m(1) \right| \ \leq \ {n_1+n_2 \over \beta^m (\beta-1) (m+1)}.$$
\endproclaim
\demo{Proof} In the case when $g(z)$ is a polynomial (that is, when $h_2(z) \equiv 1$), this is Lemma 2.2.1 of \cite{Ba16}. The proof below in the case of a rational function is very similar.

Let $\alpha_{11}, \ldots, \alpha_{1n_1}$ be the roots of $h_1$ and let $\alpha_{21}, \ldots, \alpha_{2n_2}$ be the roots of $h_2$, counting multiplicity.
Hence 
$$h_1(z)=h_1(0) \prod_{i=1}^{n_1} \left(1- {z \over \alpha_{1i}}\right) \quad \text{and} \quad h_2(z)=h_2(0) \prod_{j=1}^{n_2} \left(1- {z \over \alpha_{2j}}\right),$$
where
$$| \alpha_{1i}| \ \geq \ \beta \quad \text{for} \quad i=1, \ldots, n_1 \quad \text{and} \quad |\alpha_{2j}| \ \geq \ \beta \quad \text{for} \quad j=1, \ldots, n_2.$$
Then 
$$f(z)= f(0) +  \sum_{i=1}^{n_1} \ln \left(1-{z \over \alpha_{1i}}\right) - \sum_{j=1}^{n_2} \ln \left(1 - {z \over \alpha_{2j}}\right),$$
where we choose the branch of the logarithm so that $\ln 1=0$.

Approximating the logarithms by their Taylor polynomials, we obtain 
$$\ln \left(1- {1\over \alpha_{1i}}\right)=-\sum_{k=1}^m {1 \over k \alpha_{1i}^k} + \eta_{1i} \quad \text{and} \quad 
\ln \left(1- {1\over \alpha_{2j}}\right)=-\sum_{k=1}^m {1 \over k \alpha_{2j}^k} +\eta_{2j},$$
where 
$$\left| \eta_{1i} \right| = \left| \sum_{k=m+1}^{\infty} {1 \over k \alpha_{1i}^k} \right| \ \leq \ {1 \over m+1} \sum_{k=m+1}^{\infty}  {1 \over \beta^k} = 
{1 \over (m+1) \beta^m (\beta-1)}$$ 
for $i=1, \ldots, n_1$ and, similarly,
$$\left| \eta_{2j} \right| = \left| \sum_{k=m+1}^{\infty} {1 \over k \alpha_{2j}^k} \right| \ \leq \ {1 \over m+1} \sum_{k=m+1}^{\infty}  {1 \over \beta^k} = 
{1 \over (m+1) \beta^m (\beta-1)}$$
for $j=1, \ldots, n_2$.

Since 
$$T_m(1)= -\sum_{i=1}^{n_1} \sum_{k=1}^m {1 \over k \alpha_{1i}^k} + \sum_{j=1}^{n_2} \sum_{k=1}^m {1 \over k \alpha_{2j}^k},$$
the proof follows.
{\hfill \hfill \hfill} \qed
\enddemo

\proclaim{(2.2) Corollary} For $0< \epsilon < 1$, under the conditions of Lemma 2.1, we have 
$$\left| f(1)-T_m(1)\right| \ \leq \ \epsilon$$
provided 
$$m \ \geq \ {c \over \beta-1} \ln {n_1+n_2 \over \epsilon (\beta-1)}$$
where $c >0$ is an absolute constant. 
\endproclaim 

\demo{Proof} Follows by Lemma 2.1.
{\hfill \hfill \hfill} \qed
\enddemo

To compute the value of $T_m(1)$ in Lemma 2.1 and Corollary 2.2, we need to compute the derivatives $f^{(k)}(0)$ for $k=0, 1, \ldots, m$. This, in turn, reduces to computing the derivatives $g^{(k)}(0)$ for $k=0, 1, \ldots, m$, as explained in Section 2.2.2 of \cite{Ba16}. For completeness, we describe the procedure here.

\subhead (2.3) Computing $f^{(k)}(0)$ from $g^{(k)}(0)$ \endsubhead  We have 
$$f'(z)={g'(z) \over g(z)} \quad \text{from which} \quad g'(z) = f'(z) g(z)$$
and hence
$$g^{(k)}(0)=\sum_{j=0}^{k-1} {k-1 \choose j} f^{(k-j)}(0) g^{(j)}(0) \quad \text{for} \quad k=1, \ldots, m. \tag2.3.1$$
Now, (2.3.1) is a triangular system of linear equations in the variables $f^{(k)}(0)$ for $k=1, \ldots, m$ with diagonal coefficients $g^{(0)}(0)=g(0) \ne 0$, so the matrix of the system is invertible. Given the values of $g(0)$ and $g^{(k)}(0)$ for $k=1, \ldots, m$, one can compute the values of $f^{(k)}(0)$ for $k=1, \ldots, m$ in 
$O(m^2)$ time. This is, of course, akin to computing cumulants of a probability distribution from its moments.

Finally, we employ a rational transformation.

\proclaim{(2.4) Lemma} For real $0 < \rho < 1$, let
$$\xi=\xi_{\rho}=1-\sqrt{\rho \over 1+\rho}, \quad  \beta=\beta_{\rho}=\xi^{-1} \ \geq \ 1+\sqrt{\rho}.$$
and let
$$\psi=\psi_{\rho}(z)={\rho \over (1-\xi z)^2} - \rho, \quad \psi: \widehat{\Bbb C} \longrightarrow \widehat{\Bbb C},$$
be a rational function.
Then $\psi(0)=0$, $\psi(1)=1$ and the image of the disc
$$\beta {\Bbb D}=\left\{ z \in {\Bbb C}: \ |z| \ < \ \beta \right\}$$ 
under $\psi$ does not intersect 
the ray
$$\left\{ z \in {\Bbb C}: \ \Im \thinspace{z} =0 \quad \text{and} \quad \Re\thinspace z \ \leq \ -{3 \rho \over 4} \right\}.$$
\endproclaim
\demo{Proof} Clearly, $\psi(0)=0$ and $\psi(1)=1$. For $z \in \beta {\Bbb D}$, we have $|\xi z| < 1$ and hence 
$$\arg {1 \over 1-\xi z} \ < \ {\pi \over 2}.$$
Therefore the image of $\beta {\Bbb D}$ under the map 
$$z \longmapsto {\rho \over (1-\xi z)^2} \tag2.4.1$$ 
does not contain the non-positive real ray 
$$R_-=\left\{z\in {\Bbb C}:\  \Im \thinspace z =0, \ \Re\thinspace z \ \leq \ 0 \right\}.$$
The real values of the map (2.4.1) on the disc $\beta {\Bbb D}$ are attained when $z$ is real, and are larger than $\rho/4$, which is attained when $z=-\beta$.

The proof now follows.
{\hfill \hfill \hfill} \qed 
\enddemo

Now we are ready to prove Theorem 1.1.
\subhead (2.5) Proof of Theorem 1.1 \endsubhead Let 
$\rho = 4\delta/3$ and let 
$\psi=\psi_{\rho}: \ \widehat{\Bbb C} \longrightarrow \widehat{\Bbb C}$ be the corresponding rational transformation of Lemma 2.4.
We consider the composition 
$$g(z)=p\left(\psi(z)\right).$$ 
Clearly, $g(0)=p(0)$ and $g(1)=p(1)$.
Let 
$$\xi =1-\sqrt{\rho \over 1+\rho} \quad \text{and} \quad \beta = \xi^{-1} \ \geq \ 1+ \sqrt{4 \delta \over 3},$$
as in Lemma 2.4.
 Since the image $\psi\left(\beta {\Bbb D}\right)$ does not intersect the ray 
$$R=\left\{ z \in {\Bbb C}: \ \Im\thinspace z=0 \quad \text{and} \quad \Re\thinspace z \ \leq \ -\delta \right\},$$
we conclude that 
$$g(z) \ne 0 \quad \text{provided} \quad |z| \ < \ \beta.$$
For some polynomials $h_1(z)$ and $h_2(z)$, we can write  
$$g(z) = {h_1(z) \over h_2(z)} \quad \text{where} \quad h_2(z) = (1-\xi z)^{2n} \quad \text{and} \quad \deg h_2(z) \ \leq \ 2n.$$
Let us choose a branch of 
$$f(z)=\ln g(z) \quad \text{for} \quad |z| < \beta$$
and let $T_m(z)$ be the Taylor polynomial of $f$ of degree $m$, computed at $z=0$. From Corollary 2.2, we have 
$$\left| T_m(1) - f(1)\right| = \left| T_m(1) - \ln p(1)\right| \ \leq \ \epsilon,$$
as long as 
$$m \ \geq \ {c \over \sqrt{\delta}} \ln {n \over \epsilon  \sqrt{\delta}}$$
for some absolute constant $c > 0$. 

It remains to show how to compute the values $f^{(k)}(0)$ for $k=0, \ldots, m$ from the coefficients $a_k$, $k=0, 1, \ldots, m$, of the polynomial $p$.
First, we compute the values $g^{(k)}(0)$ for $k=0, \ldots, m$. To that effect, let 
$$p_{[m]}(z)=\sum_{k=0}^m a_k z^k$$ 
be the truncation of the polynomial $p$ and let 
$$\psi_{[m]}(z)= \rho \sum_{k=1}^m (k+1) \xi^k z^k$$
be the truncation of the Taylor series expansion of $\psi(z)$ in the disc $\beta {\Bbb D}$. Note that since $\psi(0)=0$, the constant term of the expansion is $0$.
We then compute the composition 
$$p_{[m]}\left(\psi_{[m]}(z)\right)$$ 
and discard the terms of degree higher than $m$. A fast way to do it is by Horner's method, successively computing 
$$\left( \ldots \left(\left(a_m \psi_{[m]}(z) + a_{m-1}\right) \psi_{[m]}(z)\right)+ \ldots \right) \psi_{[m]}(z)+ a_0$$
and discarding monomials of degree higher than $m$ on each step.
This gives us the Taylor polynomial of degree $m$ of $g(z)$, computed at $z=0$. We then compute the derivatives $f^{(k)}(0)$ 
as in Section 2.3.
{\hfill \hfill \hfill} \qed

To prove Theorem 1.4, we use a different (simpler) rational transformation.

\proclaim{(2.6) Lemma} For real $0 < \rho < 1$,  let
$$\xi=\xi_{\rho}={1 \over 1+\rho}, \quad  \beta=\beta_{\rho}=\xi^{-1} =1+\rho.$$
and let
$$\psi=\psi_{\rho}(z)={\rho \over 1-\xi z} - \rho, \quad \psi: \widehat{\Bbb C} \longrightarrow \widehat{\Bbb C},$$
be a rational function.
Then $\psi(0)=0$, $\psi(1)=1$ and the image of the disc
$$\beta {\Bbb D}=\left\{ z \in {\Bbb C}: \ |z| \ < \ \beta \right\}$$ 
under $\psi$ does not intersect 
the half-plane
$$\left\{ z \in {\Bbb C}: \quad \Re\thinspace z \ \leq \ -\rho \right\}.$$
\endproclaim
\demo{Proof} Clearly, $\psi(0)=0$ and $\psi(1)=1$. For $z \in \beta {\Bbb D}$, we have $|\xi z| <1$ and hence 
$$\arg {1 \over 1-\xi z} \ < \ {\pi \over 2}.$$
Therefore, the image of $\beta {\Bbb D}$ under the map 
$$z \longmapsto {1 \over 1-\xi z}$$
does not intersect the half-plane $\Re \thinspace z \leq 0$. 

The proof now follows.
{\hfill \hfill \hfill} \qed
\enddemo

\subhead (2.7) Proof of Theorem 1.4 \endsubhead We define the transformation $\psi=\psi_{\delta}$ as in Lemma 2.6, consider the composition 
$g(z)=p(\psi(z))$ and proceed as in the proof on Theorem 1.1 in Section 2.5 with straightforward modifications.
{\hfill \hfill \hfill} \qed

\subhead (2.8) Remark: approximating $p'(1)$ \endsubhead  It follows from Rolle's Theorem that if $p(x)=a_0+a_1 x+ \ldots + a_n x^n$  is a non-constant polynomial satisfying the conditions of Theorem 1.1 then 
$p'(x)=a_1 + 2 a_2 x + \ldots + n a_n x^{n-1}$ also satisfies the conditions of Theorem 1.1. Similarly, it follows from the Gauss-Lucas Theorem that if a non-constant polynomial $p(z)$ satisfies the conditions of Theorem 1.4, then so does $p'(z)$.

\subhead (2.9) Possible ramifications \endsubhead To prove Theorem 1.6, for $0 < \rho < 1$, we define 
$$\xi=\xi_{\rho}=1-e^{-{1 \over \rho}}, \quad \beta = {1-e^{-1-{1 \over \rho}} \over 1- e^{-{1\over \rho}}}  \ > \ 1 $$ and consider the map
$$\psi=\psi_{\rho}(z)=\rho \ln {1 \over 1- \xi z} \quad \text{where} \quad |z| \ < \ \beta.$$
Then $\psi(0)=0$, $\psi(1)=1$ and if $\rho=\rho(\delta) >0$ is chosen small enough, the image $\psi\left(\beta {\Bbb D}\right)$ lies in the $\delta/2$-neighborhood of 
$[0, 1] \subset {\Bbb C}$. However, $\psi$ is no longer a rational transformation. Instead of $\psi$, we use its scaled Taylor polynomial approximation $\tilde{\psi}$ such that we still have 
$\tilde{\psi}(0)=0$, $\tilde{\psi}(1)=1$ and $\tilde{\psi}$ maps the disc $\beta {\Bbb D}$ inside the $\delta$-neighborhood of $[0, 1] \subset {\Bbb C}$, see Section 2.2 of \cite{Ba16} for details. Then the proof proceeds as in Section 2.5.

Suppose now that $p(z)$ is a polynomial of degree $n$ such that $p(z) \ne 0$ whenever $z$ lies in the $\delta$-neighborhood of the sector 
$$S_{\alpha}=\left\{z \in {\Bbb C} \setminus \{0\}: \quad | \arg z| \leq \alpha \right\}$$ 
for some fixed $\alpha >0$ and some $\delta >0$. In this case, for $0 < \rho < 1$ define
$$\xi=\xi_{\rho}=1-\left({ \rho \over 1+\rho}\right)^{\pi / 2 \alpha}, \quad \beta =\beta_{\rho}=\xi^{-1}\left(1- {1 \over 2} \left({ \rho \over 1+\rho}\right)^{\pi / 2 \alpha}\right) \ > \ 1$$
and consider the (well-defined) map
$$\psi(z)=\psi_{\rho}(z)=\rho \left(1 - \xi z \right)^{-{2 \alpha /\pi}}-\rho \quad \text{for} \quad |z| < \beta.$$
We observe that $\psi(0)=0$, $\psi(1)=1$ and one can show that by choosing $\rho$ small enough, we can make sure that the image of the disc $\beta {\Bbb D}$ lies in a
prescribed neighborhood of the sector $|\arg z | \leq \alpha$. We then can use a sufficiently accurate polynomial approximation $\tilde{\psi}$ to $\psi$ 
to show that the value of $p(1)$, up to relative error $\epsilon$, is determined by the lowest
$${c \over \delta^{\pi/2\alpha}}  \ln {n \over  \epsilon \delta^{\pi/2\alpha}}$$
coefficients of $p$, where $c >0$ is an absolute constant.

\head 3. Proofs of Theorems 1.2, 1.3 and 1.5 \endhead

\subhead (3.1) Proof of Theorem 1.2 \endsubhead Given a graph $G$, we define its {\it matching polynomial} by 
$$p(x)=1 + \sum_{k=1}^{v(G)} m_k(G) x^k.$$
The Heilmann-Lieb Theorem \cite{HL72}, see also \cite{GG81}, asserts that the roots $x$ of $p(x)$ are real and satisfy $x \leq -\delta$ for 
$$\delta = {1 \over 4\left(\Delta-1\right)} \quad \text{for} \quad \Delta=\max\left\{ \Delta(G), 2 \right\}.$$
The proof now follows from Theorem 1.1. 
{\hfill \hfill \hfill} \qed

\subhead (3.2) Proof of Theorem 1.3 \endsubhead Given a graph $G$, we define its {\it independence polynomial} by 
$$p(x)=1 + \sum_{k=1}^{v(G)} i_k(G) x^k.$$ Chudnovsky and Seymour proved \cite{CS07} that if $G$ is claw-free, then the roots of $p(x)$ are necessarily non-positive real.
On the other hand, the Dobrushin-Shearer bound, cf. \cite{SS05}, states that the roots $z$ of the independence polynomial of any graph $G$ satisfy 
$$|z| \ \geq \ {({\Delta -1)}^{\Delta-1} \over {\Delta}^{\Delta}} ={1 \over \Delta e} \left(1 + O\left({1 \over \Delta}\right)\right) \quad \text{as} \quad \Delta \longrightarrow \infty,$$
where $\Delta=\max\{2 , \Delta(G)\}$. The proof now follows from Theorem 1.1.
{\hfill \hfill \hfill} \qed

\subhead (3.3) Proof of Theorem 1.5 \endsubhead Given a graph $G$, we define its {\it unbranched subgraph polynomial} by 
$$p(z)=1+\sum_{k=1}^{v(G)} u_k(G) z^k.$$ 
Ruelle proved \cite{R99a}, \cite{R99b}, see also \cite{Wa09}, that all roots $z$ of $p(z)$ satisfy 
$$ \Re\thinspace z \ \leq \ -{2 \over \Delta(\Delta-1)^2} \quad \text{for} \quad \Delta=\max\left\{2, \Delta(G)\right\}.$$
The proof now follows from Theorem 1.4.
{\hfill \hfill \hfill} \qed

\subhead (3.4) Estimating averages \endsubhead We note that the value 
${p'(1) \over  p(1)}$
is interpreted the average number of edges in a matching in Theorem 1.2, the average number of vertices in an independent set in Theorem 1.3 and the average number of edges in a unbranched subgraph in Theorem 1.5. It follows from Remark 2.8 that we can estimate the averages within relative error $\epsilon >0$ by inspecting the matchings, 
independent sets and unbranched subgraphs of pretty much the same size as prescribed by Theorems 1.2, 1.3 and 1.5, though with different absolute constants.
Similarly, by computing ${p''(1) \over p(1)}$ we can estimate the second moment, etc.

\subhead (3.5) Possible ramifications \endsubhead For each vertex $w$ of a graph $G$, let us choose a set $A_w$ of allowable degrees of subgraphs. 
Wagner proved \cite{Wa09} that if $0 \in A_w \subset \{0, 1, 2 \}$ for all vertices $w$, then the corresponding subgraph counting polynomial is non-zero in 
the sector 
$$S_{\pi/3}=\left\{ z \in {\Bbb C}:\ |\arg z| < {\pi \over 3}\right\}$$ and is also non-zero in a $\delta$-neighborhood of $z=0$ for some $\delta=\Omega(1/\Delta(G))$.
Using the approach sketched in Section 2.8, one can show that within relative error $\epsilon >0$, the total number of such subgraphs is determined by the numbers of 
subgraphs with 
$$k=\left(\Delta(G)\right)^{O(1)} \ln {v(G) \over \epsilon}$$
edges.

\enddocument
\end